\documentclass[preprint,11pt]{amsart}
\usepackage{amscd,amssymb,amsopn,amsmath,amsthm,graphics,amsfonts,enumerate,verbatim,calc}
\usepackage[dvips]{graphicx}
\usepackage[colorlinks=true,linkcolor=red,citecolor=blue]{hyperref}
\newcommand{\ra}{\rightarrow}

\newcommand{\expx}{\exp_{x}}
\newcommand{\expxv}{\exp_{x}^{-1}}
\newcommand{\expxy}{\exp_{x}^{-1}{y}}
\newcommand{\expyx}{\exp_{y}^{-1}{x}}
\newcommand{\norm}{\|}

\newcommand{\rbar}{\bar{r}}
\newcommand{\xbar}{\bar{x}}

\newcommand{\hess}{{\rm{Hess }}\;}
\newcommand{\unp}{{\rm{Unp }}}
\newcommand{\reach}{{\rm{reach }}}

\newcommand{\cone}{{\rm{cone}}}
\newcommand{\spa}{{\rm{span}}}

\newcommand{\li}{\langle}%%$left inner product
\newcommand{\ri}{\rangle}%%%right inner product

\newcommand{\RR}{\mathbb{R}}

\newtheorem{theorem}{Theorem}[section]
\newtheorem{corollary}[theorem]{Corollary}
\newtheorem{lemma}[theorem]{Lemma}

\newtheorem{proposition}[theorem]{Proposition}

\newtheorem{definition}[theorem]{Definition}
\newtheorem{example}[theorem]{Example}

\theoremstyle{plain}

\newcommand{\re}[1]{{\color{black}{#1}}}
\newcommand{\blu}[1]{{\color{black}{#1}}}

\numberwithin{equation}{section}
\begin{document}
\title[Minimizing curves in prox-regular sets]{Minimizing curves in prox-regular subsets of Riemannian manifolds}
%\author[M. R. Pouryayevali and H. Radmanesh]{M. R. Pouryayevali$^{\dag, *}$, and H. Radmanesh$^{\dag}$,}
\address{Department of Pure Mathematics, Faculty of Mathematics and Statistics, University of Isfahan, Isfahan, 81746-73441, Iran}
\email{pourya@math.ui.ac.ir} \email{h.radmanesh@sci.ui.ac.ir}
%\thanks{$*$Corresponding author.}

%\subjclass[2010]{58C20, 58C06, 49J52}
\author[M. R. Pouryayevali and H. Radmanesh]{Mohamad R. Pouryayevali and Hajar Radmanesh}
%\thanks{\emph{2010 Mathematics Subject Classification: 58C20, 58C06, 49J52}\\ $\dag$ Department of Mathematics, University of Isfahan, P.O. Box 81745-163, Isfahan, Iran\\
%Emails: pourya@math.ui.ac.ir, h.radmanesh@sci.ui.ac.ir\\
%$*$Corresponding author.}

\maketitle

%%%%%%%%%%%%%%%%%%%%%%%%%%%%%%%%%%%%%%%%%%%%%%%%%%%%%%%%%%%%%%%%%%%%%%%%%%%%%%%%%%%%%%%%%%%%%%%%%%%%%%%%%%%%%%%%%%%%%%%%%%%%%%%%%
%%%%%%%%%%%%%%%%%%%%%%%%%%%%%%%%%%%%%%%%%%%%%%%%%%%%%%%%%%%%%%%%%%%%%%%%%%%%%%%%%%%%%%%%%%%%%%%%%%%%%%%%%%%%%%%%%%%%%%%%%%%%%%%%%
\begin{abstract}
We obtain a characterization of the proximal normal cone to a prox-regular subset of a Riemannian manifold. Moreover,  some properties of Bouligand tangent cones to prox-regular sets are described.
 We prove that for a prox-regular subset $S$ of a Riemannian manifold,  the metric projection $P_S$ to $S$ is locally Lipschitz on an open neighborhood of $S$ and it is directionally differentiable at boundary points of $S$. Finally, a  necessary condition for a curve to be a minimizing curve in a prox-regular set is derived.

\end{abstract}
%%%%%%%%%%%%%%%%%%%%%%%%%%%%%%%%%%%%%%%%%%%%%%%%%%%%%%%%%%%%%%%%%%%%%%%%%%%%%%%%%%%%%%%%%%%%%%%%%%%%%%%%%%%%%%%%%%%%%%%%%%%%%%%%
%%%%%%%%%%%%%%%%%%%%%%%%%%%%%%%%%%%%%%%%%%%%%%%%%%%%%%%%%%%%%%%%%%%%%%%%%%%%%%%%%%%%%%%%%%%%%%%%%%%%%%%%%%%%%%%%%%%%%%%%%%%%%%%%
\section{Introduction}\label{sec1}
%\section{Preliminaries and notations}\label{sec2}
Closed subsets of  Hilbert spaces satisfying an external sphere condition with uniform radius have been studied as generalizations of convex sets, mostly in relation to uniqueness of the metric projection and smoothness of the distance function. In the fundamental paper \cite{Federer} where the finite dimensional case is considered, these sets were called sets with positive reach. Then various \blu{equivalent} definitions related to  this property have been presented independently by several authors; see  \cite{colombo1,RO} and the references therein.  Among them, one can mention \blu{the notions of $\varphi$-convexity (as titled $p$-convexity)
 %$p$-convex sets
 and prox-regularity of sets which were introduced in \cite{colombo2} and \cite{RO}, respectively.}  It was shown in \cite{canino} that certain properties which hold globally for convex sets are still valid locally for \blu{$\varphi$-convex sets}.

Differentiability properties of the metric projection onto closed convex sets are of interest  in sensitivity analysis of  variational inequalities and optimal control problems. Moreover, the regularity of the metric projection  onto a sufficiently regular submanifold $M$ of $\mathbb{R}^{n}$ as well as the regularity of the corresponding distance
function have significant role in various aspects of analysis; see \cite{Leobacher}. A  classical example is the Dirichlet problem for quasilinear partial differential equations, where the manifold of interest is the boundary of the underlying domain; see, for instance, \cite{Gilberg}.

The example presented  by J. Kruskal \cite{kruskal} shows that, in general such a projection is not directionally differentiable, even in finite dimensional spaces. By directionally differentiable at a point we mean that the directional derivative
exists for all directions through that point. This is weaker than the existence of the
gradient at that point.

The problem of differentiability of the metric
projection  for a closed locally convex subset $S$ of a finite dimensional Riemannian
manifold $M$ was studied in \cite{Walter} and  it was proved in \cite{Greene} that for a closed totally convex subset $S\subset M$,   there exists an open set $W$ containing
$S$ such that the metric projection is locally Lipschitz on $W$.

 In \cite{Barani} the notion of \blu{$\varphi$-convex sets was extended  to Hadamard manifolds} and it was shown that if $S$ is a $\varphi$-convex subset of an infinite-dimensional Hadamard manifold $M$, then there exists a neighborhood $U$ of $S$ in $M$ such that the metric projection $P_S: U\to S$ is single-valued and locally Lipschitz. Moreover, it was proved that under the same assumptions on $S$ and $M$, there exists a neighborhood $U$ of $S$ in $M$ such that $d^{2}_S$ is $C^1$ with locally Lipschitz gradient on $U\setminus S$. \blu{On the other hand, in \cite{Hosseini on}  the notion
of prox-regular sets  was introduced on Riemannian manifolds as a subclass of regular sets. In \cite{convex} we proved that   the two classes of $\varphi$-convex sets and prox-regular sets coincide in the setting of Riemannian manifolds.}

\blu{The problem of existence and uniqueness of geodesics on a Riemannian manifold without boundary is a classical subject of differential geometry and global nonlinear analysis} and is particularly fit to a treatment by variational methods. However, in the case of Riemannian manifolds with boundary or certain subsets of a manifold without boundary,  strong irregularities appear  in the energy functional and  new techniques are needed for dealing with these problems. In \cite{canino1,canino} \blu{$\varphi$-convex subsets of a Real Hilbert space were} considered and using an infinitesimal definition of geodesics in the framework of Sobolev spaces the author characterized these geodesics as critical points of an energy functional on a suitable path spaces.

The aim of this paper is to study minimizing curves in a prox-regular subset $S$ of a Riemannian manifold $M$. To this end, we use \blu{some} powerful tools from nonsmooth analysis and an adapted variational technique. Applying the first variation formula, we give a necessary condition for an admissible  curve $\gamma:[a,b]\ra M$ in $S$ to be minimizing. Indeed, this curve \blu{has} the property that %almost all its covariant derivatives are contained in the proximal normal cones to $S$ at
\[
D_t\dot{\gamma}(t)\in N^P_S\left(\gamma(t)\right),
 \]
for every $t\in [a,b]$ except for finitely many points, provided that $S$ has a $C^2$ boundary, where $N^P_S(x)$ is the proximal normal cone at $x\in S$.  To prove this result, we  address the problem of the directional differentiability of the metric projection $P_S$ at boundary points of $S$. Employing Shapiro's variational principle \cite{Shapiro1}, we show that for a prox-regular subset $S$ of a Riemannian manifold $M$, the projection map $P_S$ is locally Lipschitz on an open neighborhood of $S$ which generalizes the result of \cite{Barani} to the Riemannian setting. Moreover, we prove that $P_S$ is directionally differentiable at boundary points of $S$. We  also obtain a useful characterization of Bouligand tangent cone to a prox-regular set.

The paper is organized as follows. In Section \ref{sec2}  we present some basic constructions and
preliminaries in Riemannian geometry and nonsmooth analysis, widely used in the sequel.  Section \ref{sec3} is devoted to the study of Bouligand and proximal normal cones. %Then we prove that every prox-regular set is locally connected. This result is used to
\blu{Then we} obtain a characterization of the proximal normal cone to a prox-regular set.
%In Section \ref{sec4}  certain topological properties of $\varphi$-convex sets are discussed.
We also show that $P_S$ is a locally Lipschitz retraction from a neighborhood of $S$ to $S$. In Section \ref{sec5} differentiability properties of the metric projection $P_S$ to a prox-regular subset $S$ of a Riemannian manifold are investigated  which leads to a characterization of Bouligand tangent cone. Section \ref{sec6} is concerned with the necessary condition for a curve $\gamma$ to be a minimizing curve in a prox-regular set whose boundary is a $C^2$ submanifold of $M$. Moreover,  some relevant examples are presented.

%%%%%%%%%%%%%%%%%%%%%%%%%%%%%%%%%%%%%%%%%%%%%%%%%%%%%%%%%%%%%%%%%%%%%%%%%%%%%%%%%%%%%%%%%%%%%%%%%%%%%%%%%%%%%%%%%%
%%%%%%%%%%%%%%%%%%%%%%%%%%%%%%%%%%%%%%%%%%%%%%%%%%%%%%%%%%%%%%%%%%%%%%%%%%%%%%%%%%%%%%%%%%%%%%%%%%%%%%%%%%%%%%%%%%
 \section{Preliminaries and notations}\label{sec2}
Let us recall some notions of Riemannian manifolds and nonsmooth analysis; see, e.g.,
\cite{Clark book a,Do Carmo,sakai}. Throughout this paper, $(M,g)$ is a
finite-dimensional Riemannian manifold endowed with a Riemannian metric $g_x=\li .,.\ri_x$ on each tangent space $T_xM$ and  $\nabla$ is the
Riemannian connection of $g$.  For every $x, y \in M$, the
Riemannian distance from $x$ to $y$ is denoted by $d(x,y)$. Moreover, $B(x,r)$ and $\overline{B}(x,r)$ signify the
open and closed metric ball centered at $x$ with radius $r$, respectively. For a smooth curve $\gamma : I \rightarrow M$ and $t_0,t\in I$, the notation  $L^{\gamma}_{t_0t}$ is used for the parallel transport along $\gamma$ from $\gamma(t_0)$ to $\gamma(t)$. When $\gamma$ is the unique minimizing geodesic joining $\gamma(t_0)$ to $\gamma(t)$, we use $L_{t_0t}$ instead of $L^{\gamma}_{t_0t}$. Furthermore for a smooth vector field $X$ along $\gamma$, $D_t X$ is the covariant derivative of $X$ along $\gamma$.

For $x\in M$, let $r(x)$ be the convexity radius at $x$. Then the function $x\mapsto r(x)$ from $M$ to $\mathbb{R}^{+} \cup \{
+\infty \}$ is continuous; see \cite{sakai}.  The map $\exp_x : U_x \to M$ will stand for the exponential map
at $x$, where $U_x $ is an open subset of the tangent space $T_xM$
containing $0_{x} \in T_xM$.
 Note that if $x$ and $y$ belong to a convex set, then both
$\expxy$ and $\expyx$ are defined and
  $$\norm \expxy \norm = d(x,y) = \norm \expyx \norm.$$
Moreover,
\[
 L_{xy}\left(\expxy\right) = - \expyx.
 \]
For a fixed point $z\in M$, the function $\phi:M\ra \RR$ defined by $\phi(x)=d^2(x,z)$ is $C^{\infty}$ on any convex neighborhood  of $z$ and for every $x$ in a convex neighborhood of $z$,  $\nabla \phi(x)= - 2\exp_x^{-1}z$.

Let $S$ be a nonempty closed subset of  $M$. The proximal normal cone to $S$ at $x\in S$,
is denoted by $N^{P}_{S}(x)$ and $\xi\in N^{P}_{S}(x)$ if and only if there exists  $\sigma
>0$ such that
\[
\li \xi, \expxy \ri \leq\sigma\;d^2(x,y),
\]
for every $y\in U\cap S$, where $U$ is a convex neighborhood of $x$.
The metric projection to $S$, denoted by $P_S$, is defined by
\[
P_S(z)= \left\{ x\in S : d(x,z)=\inf_{y\in S} \, d(y,z)\right \}\quad \forall z\in M.
\]
Moreover, $\unp(S)$ is considered as the set of all points $z\in M$ with the property that $P_S(z)$ is single-valued.
Then according to \cite[Lemma 4.11]{convex}, the projection map $P_S: \unp(S)\rightarrow S$ is continuous.  For every $x\in S$ we also define
\[
\reach(S,x):=\sup\left\{r\geq 0: B(x,r)\subseteq \unp(S)\right\},
\]
 It is worth mentioning that the function $x\mapsto \reach(S,x)$ is continuous on $S$; see \cite{EFP,Klein2} for
more details.

In order to deduce the  Lipschitz property and directional differentiability of $P_S$, we use the following variational principal by A. Shapiro \cite{Shapiro1}. Let $f,g:X\ra \mathbb{R}$ be two functions on a Hilbert space $X$ and $S,T\subset X$. Consider the optimization problems
\begin{equation}\label{c}
  \min_{x\in S}f(x)
\end{equation}
and
\begin{equation}\label{cc}
  \min_{x\in T}g(x).
\end{equation}
Let $x_0$ and $\xbar$ be some optimal solutions of (\ref{c}) and (\ref{cc}), respectively and suppose that there exist a neighborhood $W$ of $x_0$ and  $\alpha>0$ such that for every $x\in S\cap W$,
\begin{equation}\label{c1}
  f(x)\geq f(x_0)+\alpha \norm x-x_0\norm ^2.
\end{equation}
Also, suppose that $\xbar\in W$ and $f$ and $g$ are Lipschitz on $W$ with Lipschitz constants $k_1$ and $k_2$, respectively. Then
\begin{equation}\label{c2}
  \norm \xbar-x_0\norm\leq \alpha^{-1}\kappa+2\delta_1+\alpha^{-1/2}\left(k_1\delta_1+k_2\delta_2\right)^{1/2},
\end{equation}
where $\kappa$ is a Lipschitz constant of  $h(x)=g(x)-f(x)$ on $W$ and
\[
\delta_1=\sup_{x\in T\cap W}d(x,S\cap W),
\]
\[
\delta_2=d(x_0,T\cap W).
\]
%%%%%%%%%%%%%%%%%%%%%%%%%%%%%%%%%%%%%%%%%%%%%%%%%%%%%%%%%%%%%%%%%%%%%%%%%%%%%%%%%%%%%%%%%%%%%%%%%%%%%%%%%%%%%%%%%%%%%%%%%%%%%%%%%
%%%%%%%%%%%%%%%%%%%%%%%%%%%%%%%%%%%%%%%%%%%%%%%%%%%%%%%%%%%%%%%%%%%%%%%%%%%%%%%%%%%%%%%%%%%%%%%%%%%%%%%%%%%%%%%%%%%%%%%%%%%%%%%%%

\section{Local Lipschitzness of metric projection}\label{sec3}
In this section we first derive  some properties of Bouligand tangent cones to prox-regular sets which we need in the sequel. Let us begin by recalling some required definitions; see \cite{Hosseini on,convex}.

\blu{The closed subset $S$ of $M$ is said to be prox-regular at $\bar{x}\in S$ if there exist $\varepsilon>0$ and $\sigma>0$ such that $B(\bar {x}, \varepsilon)$ is convex and for every  $x\in S\cap B(\bar{x}, \varepsilon)$ and $v\in N^P_S(x)$ with $\norm v\norm <\epsilon$,
  \[
    \li v, \expxy \ri \leq\sigma\;d^2(x,y) ~ \textrm{for every}\ y\in S\cap B(\bar{x}, \epsilon).
  \]
Moreover, $S$ is called prox-regular if it is prox-regular at each point of $S$.

In \cite[Theorem 3.4]{convex}, we proved that   every $\varphi$-convex subset of a Riemannian manifold $M$ is prox-regular and conversely,  for every prox-regular subset $S$ of $M$ there exists a continuous function $\varphi : S\rightarrow [0,\infty)$ such that $S$ is  $\varphi$-convex.  %Therefore whenever the function $\varphi$ plays no role explicitly, the expression ``$\varphi$-convex'' may be replaced by ``prox-regular''.
Recall that a closed subset $S\subset M$ is called $\varphi$-convex if for every $x\in S$ and $v\in N^P_S(x)$
  \[
   \li v,\exp_x^{-1}y \ri \leq \varphi(x)\norm v\norm d^2(x,y),
   \]
for every $y\in U\cap S$, where $U$ is a convex neighborhood of $x$ and $\varphi : S\rightarrow [0,\infty)$ is a continuous function. Note that this definition is independent of the choice of any convex neighborhood of $x$.}

Let $S\subset M$ be a closed subset and $x \in S$. The Bouligand (or contingent) tangent cone to $S$ at $x$ is defined as
\[
T^B_S(x):=\left\{\lim_{i\ra \infty}\frac{\exp_x^{-1} z_i}{t_i} : z_i\in U\cap S, z_i\ra x\  \textrm{and}\
t_i\downarrow 0\right\},
\]
where $U$ is a convex neighborhood of $x$ in $M$. It was shown in \cite{Hosseini on} that when $S$ is prox-regular,   $T^B_S(x)$ is a convex cone for every $x\in S$.

\begin{lemma}\label{pol}
Let $S\subseteq M$ be a prox-regular set and $x\in S$. Then
\begin{itemize}
  \item[(i)] $T^B_S(x)=\left(N^P_S(x)\right)^{\circ}$,
  \item[(ii)] $\left(T^B_S(x)\right)^{\circ}=N^P_S(x)$.
\end{itemize}
\end{lemma}
\begin{proof}
  Assertion (i) can be obtained from \cite[Lemma 3.7]{Hosseini on}. Indeed, we have
  \[
  T^C_S(x)\subseteq T^B_S(x)\subseteq \left(N^P_S(x)\right)^{\circ}=\left(N^C_S(x)\right)^{\circ}=T^C_S(x),
  \]
  where $T^C_S(x)$ and $N^C_S(x)$ are (Clarke) tangent and normal cone to $S$ at $x$, respectively.

  Assertion (ii) follows from the fact that $N^P_S(x)$ is closed and convex. Hence $\left(\left(N^P_S(x)\right)^{\circ}\right)^{\circ}= N^P_S(x)$. \qed
\end{proof}

According to \cite[Proposition 4.2]{convex}, for every point $x$ in a closed prox-regular subset $S$ of $M$, $\reach(S,x)>0$. This property of prox-regular sets helps us to prove the following topological property of these sets.

\begin{lemma}\label{lconnected}
  If $S$ is a closed %prox-regular sub
  \blu{set with the property that $\reach(S,x)>0$ for every $x\in S$}, then $S$ is locally connected.
\end{lemma}
\begin{proof}
  Let $x\in S$ and $U$ be an open neighborhood of $x$ in $M$. We \re{are going to} verify that there exists a neighborhood $V$ of $x$ in $M$ such that $V\subseteq U$ and $V\cap S$ is connected.

  If this fails to be the case, then for all positive integer $n$ large enough so that $B\left(x,1/n\right)$ is convex and $B\left(x,1/n\right)\subseteq U$, the set $S_n:=S\cap B\left(x,1/n\right)$ is not connected. Suppose that $A_n$ is the connected component of $S_n$ contains $x$, the set $B_n$ is another connected component of $S_n$ and $y_n$ is an arbitrary point of $B_n$. Let $\gamma: [0,1]\ra M$ be the unique minimizing geodesic joining $x,y_n$ and hence its image is entirely in  $B\left(x,1/n\right)$.

  Note that $P_S\left(\gamma(t)\right)\in S\cap B\left(x,1/n\right)$ for every $t\in [0,1]$, since
  \[
  \begin{array}{ll}
    d\left(P_S\left(\gamma(t)\right),x\right) & \leq d\left(P_S\left(\gamma(t)\right),\gamma(t)\right)+d\left(\gamma(t),x\right) \\
     & \leq d\left(y_n,\gamma(t)\right)+d\left(\gamma(t),x\right) \\
     & =d\left(y_n,x\right)<1/n.
  \end{array}
\]
  We now claim that the image of $\gamma$ on $[0,1]$ is not entirely in $\unp(S)$. Otherwise, the continuity of $P_S$ on $\unp(S)$ (\cite[Lemma 4.11]{convex}) implies that the set $P_S\left(\gamma([0,1])\right)$ is connected. Since  $P_S\left(\gamma([0,1])\right)\subseteq S_n$ and contains $x$, \re{we have} $P_S\left(\gamma([0,1])\right)\subseteq A_n$. It follows that $y_n\in A_n$ which contradicts our choice of $y_n$. Then there exists a sequence $\{z_n\}$ such that $z_n\notin \unp(S)$ and $d(x,z_n)<1/n$. It implies that $\reach(S,x)=0$ and this contradiction completes the proof of the lemma. \qed
\end{proof}

%It is worth noting that according to the proof of Lemma \ref{lconnected}, every closed set $S$ with the property that $\reach(S,x)>0$ for every $x\in S$ is locally connected.
\blu{Lemma \ref{lconnected} implies that every closed prox-regular subset of $M$ is locally connected.}

\begin{example}
A well known example of a connected set which is not locally connected is the comb space,
\[
C=\left([0,1]\times 0\right)\cup\left(K\times[0,1]\right)\cup\left(0\times [0,1]\right),
\]
in $\mathbb{R}^2$ where $K=\left\{1/n:n\in \mathbb{N}\right\}$. Note that this set is not prox-regular, because for every $x\in \left(0\times [0,1]\right)$, $\reach(C,x)=0$.
\end{example}

In the following theorem, we obtain a characterization of proximal normal cones to prox-regular subsets of $M$.

\begin{theorem}\label{C2}
  Suppose that $S$ is a closed subset of $M$ with the property that its boundary, denoted by $\partial S$, is an embedded $k$-dimensional  submanifold of $M$ and $x\in \partial S$. Then\\
  \emph{(a)} If $\partial S$ is $C^1$, then $N^P_S(x)\subseteq T^{\bot}_{x}\partial S$ where $T^{\bot}_{x}\partial S$ is the normal space to $\partial S$ at $x$.\\
   \emph{(b)} If in addition $S$ is prox-regular with nonempty interior and $\partial S$ is $C^2$, then
    %\begin{itemize}
     % \item [(i)] $N^P_{\partial S}(x)=T^{\bot}_{x}\partial S$.
      %\item [(ii)]
      there exist a neighborhood $U$ of $x$ in $M$ and a $C^2$ submersion $\psi: U\ra \mathbb{R}$ such that $U\cap \partial S=\psi^{-1}(0)$ and proximal normal cone to $S$ at $x$ \blu{is one of the following}
      \[
         \begin{array}{ll}
         N^P_S(x)&= \cone\left\{\nabla\psi(x)\right\},\\
         &\ {\textrm{or}}\\
          N^P_S(x)&= \spa \left\{\nabla\psi(x)\right\}.
         \end{array}
     \]

         %\begin{align*}
          %\hspace{1.2cm}  N^P_S(x) & =\left\{\Sigma_{i=1}^{j}\lambda_i\nabla\psi_i(x)- \Sigma_{i=j+1}^{n-k}\lambda_i\nabla\psi_i(x): \lambda_i\geq 0, \forall i%=1,\cdots,n-k\right\}\\
            % & =\cone\left\{\nabla\psi_1(x),\cdots, \nabla\psi_j(x),-\nabla\psi_{j+1}(x), \cdots, -\nabla\psi_{n-k}(x)\right\}.
         % \end{align*}
      %\item [(iii)] If $\hat{S}:=S^c\cup \partial S$, then $T^B_S(x)\cap  T^B_{\hat{S}}(x)=T_x\partial S$.
   % \end{itemize}
  \end{theorem}
  \begin{proof}
 Since $\partial S$ is an embedded $k$-dimensional  submanifold of $M$, there exists a neighborhood $U$ of $x$ in $M$ such that $U\cap \partial S$ is a level set of a   submersion $\psi:U\ra \mathbb{R}^{n-k}$, $\psi=\left(\psi_1,\cdots,\psi_{n-k}\right)$. If $\partial S$ is $C^1$, then \re{along the same lines as} the proof of \cite[Proposition 1.9]{Clark book a}, we have
 \[
 N^P_S(x)\subseteq N^P_{\partial S}(x)\subseteq \spa \left\{\nabla \psi_i(x): i=1,\cdots,n-k\right\}= T^{\bot}_{x}\partial S.
 \]

 If in addition $S$ is prox-regular with nonempty interior and $\partial S$ is $C^2$, then $\partial S$ is a \re{codimension 1} submanifold of $M$. Moreover, by Lemma \ref{lconnected}, $S$ is locally connected and hence by shrinking $U$ if necessary, we may assume that $U$ is convex and $U\cap S$ is connected. If $U\cap S^{\circ}=\emptyset$  where $S^{\circ}$ denotes the interior of $S$ (or there exists a neighborhood $V\subseteq U$ of $x$ such that $V\cap S^{\circ}=\emptyset$), then $U\cap S= U\cap \partial S$  and by \cite[Proposition 1.9]{Clark book a}, we have
 \[
 N^P_S(x)= N^P_{\partial S}(x)= \spa \left\{\nabla \psi(x)\right\}.
 \]

 Now let $U\cap S^{\circ}$ be nonempty. Since $U\cap S^{\circ}$ is connected and $U\cap \partial S=\psi^{-1}(0)$, we have
 \[
 \psi\left(U\cap S^{\circ}\right)\subseteq (-\infty,0)\ \ {\textrm{or}}\ \ \psi\left(U\cap S^{\circ}\right)\subseteq (0,+\infty).
 \]
  Replacing $\psi$ by $-\psi$ if necessary, we can assume that $\psi(y)\leq 0$ for every $y\in U\cap S$.  Let $\xi:=\lambda \nabla\psi(x)$ for some $\lambda\geq 0$. For given $\sigma>0$, we  define
 \[
 h(y):=\left\li -\xi, \exp_x^{-1}y\right\ri + \sigma d^2(x,y)+\lambda \psi(y),
 \]
  for every $y\in U$. Then $\nabla h(x)=0$ and for $\sigma$ sufficiently large, $\hess h(x)$ is positive definite because for every $v\in T_xM$ we have
  \[
  \re{\begin{array}{ll}
   \hess h(x)(v)^2 & = \frac{d^2}{dt^2}\mid_{t=0}\left(h\left(\exp_x(tv)\right)\right) \\
     & =\frac{d^2}{dt^2}\mid_{t=0}\left(\li -\xi,tv\ri+\sigma t^2\norm v\norm^2+\lambda \psi \left(\exp_x(tv)\right)\right) \\
     & =2\sigma \norm v\norm^2+ \lambda \hess \psi(x)(v)^2.
  \end{array}}
  \]
  Therefore $h$ has a local minimum at $x$ and so there exists a neighborhood $V$ of $x$ such that $V\subseteq U$ and for every $y\in V\cap S$ we have
 % \begin{align*}
 \[
  \left\li \xi, \exp_x^{-1}y\right\ri\leq \sigma d^2(x,y)+\lambda \psi(y)\leq \sigma d^2(x,y).
  \]
  It follows that $\xi \in N^P_S(x)$ which completes the proof of the theorem. \qed
  %\end{align*}
 \end{proof}

 It is worth mentioning that in part (b) of Theorem \ref{C2} if the interior of $S$ is empty, then $S=\partial S$ and by \cite[Proposition 1.9]{Clark book a} we have
  \[
 N^P_S(x)= \spa \left\{\nabla \psi_i(x): i=1,\cdots,n-k\right\}= T^{\bot}_{x}\partial S.
 \]

\begin{example}
Let $S$ be the set $\left(\{0\}\cup [1,+\infty)\right)\times  \mathbb{R}$ in $\mathbb{R}^2$ and consider the points $(0,0), (1,0)\in \partial S$. Then $S$ is prox-regular and has a smooth boundary. At the point $(0,0)$ we have $\psi(x,y)=x$, $ N^P_S(0,0)=\spa \{(1,0)\}$ and at the point $(1,0)$,
\[
\psi(x,y)=1-x \ \ {\textrm{and}} \ \ N^P_S(1,0)=\cone \{(-1,0)\}.
\]
\end{example}

In what follows, the closed set $S^c\cup \partial S$ is denoted by $\hat{S}$. Note that $\partial \hat{S}\subseteq \partial S$ and if the point $x\in \partial S$ is  such that $x\notin \partial \hat{S}$, then $x$ is the interior point of $\hat{S}$.
\begin{theorem}\label{tbsh}
   Suppose that $S$ is prox-regular and $\partial S$ is a $C^2$ submanifold of $M$. If $x\in \partial S$, then
   \begin{equation}\label{tbs}
      T^B_S(x)\cap T^B_{\hat{S}}(x)=T_x\partial S.
   \end{equation}
 \end{theorem}
   \begin{proof}
     In the case when $S^{\circ}=\emptyset$, we have $T^B_S(x)=T_x\partial S$ and $T^B_{\hat{S}}(x)=T_xM$. So we assume that the interior of $S$ is nonempty. Let $U$ and the submersion $\psi: U\ra \mathbb{R}$ be the ones applied in the proof of Theorem \ref{C2}.  If $U\cap S^{\circ}=\emptyset$ (or there exists a neighborhood $V\subseteq U$ of $x$ such that $V\cap S^{\circ}=\emptyset$), then $U\subseteq \hat{S}$ and $T^B_{\hat{S}}(x)=T_xM$. Hence the  expression (\ref{tbs}) holds.

     We now consider the case in which $U\cap S^{\circ}$ is nonempty and for every neighborhood $V$ of $x$ contained in $U$, $V\cap S^{\circ}\neq \emptyset$. Then $U\cap \partial \hat{S}=U\cap \partial S$ and we claim  that
     \[
     N^P_{\hat{S}}(x)=\cone \left\{-\nabla \psi(x)\right\}.
     \]
      Indeed, Since $\unp \left(\partial S\right)\subseteq \unp \left(\hat{S}\right)$ and $\partial S$ is a $C^2$ submanifold of $M$, for every $z\in \partial \hat{S}\subseteq \partial S$ we have
      \[
      \reach (\hat{S},z)\geq \reach(\partial S,z)>0.
      \]
      Then $\reach(\hat{S},z)>0$ for every $z\in \hat{S}$ and by Lemma \ref{lconnected}, $\hat{S}$ is locally connected. Without loss of generality, we assume that $U\cap \hat{S}$ is connected.  By the choice of $\psi$ we have $\psi(y)\leq 0$ for every $y\in U\cap S$. On the other hand, $\psi$ is a submersion on $U$ and $U\cap \partial \hat{S}=U\cap \partial S=\psi^{-1}(0)$. Then $\psi(y)\geq 0$ for every $y\in U\cap \hat{S}$ and  the claim is proved by a procedure similar to the proof of Theorem \ref{C2}. So we have
      \[
      N^P_S(x)\cup N^P_{\hat{S}}(x)=T^{\bot}_x\partial S\ \ {\textrm{and}}\ \  N^P_S(x)\cap N^P_{\hat{S}}(x)=\{0\}.
      \]

       Let us now prove the \re{expression} (\ref{tbs}). Since $U\cap \partial \hat{S}=U\cap \partial S$, the set $\hat{S}$ is prox-regular and applying Lemma \ref{pol}, we deduce that $T_x\partial S\subseteq T^B_S(x)\cap T^B_{\hat{S}}(x)$. Let $v\in T^B_S(x)\cap T^B_{\hat{S}}(x)$ and $w\in T^{\bot}_x\partial S$ be arbitrary. Without loss of generality, we assume that $w\in N^P_S(x)$. Thus $-w\in N^P_{\hat{S}}(x)$ and applying Lemma \ref{pol}, we have $\li v,w\ri\leq 0$ and $\li v,-w\ri\leq 0$. It follows that $v\in T_x\partial S$. \qed
   \end{proof}

   \begin{theorem}\label{o}
     Let $S$ be a $\varphi$-convex subset of $M$, $x\in S$ and \re{$U$ be a convex neighborhood of $x$.} Then
 \[
 d\left(\exp_{x}^{-1}y, T^B_S(x)\right)\leq \varphi(x)d^2\left(x,y\right),
 \]
 for every $y\in U\cap S$.
   \end{theorem}
 \begin{proof}
   \re {Since $T^B_S(x)$ is a closed convex subset of $T_xM$, for $y\in U\cap S$ there exists a unique vector  $v\in T^B_S(x)$ such that $d\left(\exp_{x}^{-1}y, T^B_S(x)\right)=\norm \exp_{x}^{-1}y-v\norm$.} Therefore we have
   \[
   \exp_{x}^{-1}y-v\in N^P_{T^B_S(x)}(v).
   \]
    Let us now show that $N^P_{T^B_S(x)}(v)\subseteq N^P_S(x)$. Clearly, $\li \xi,v\ri=0$ for every $\xi \in N^P_{T^B_S(x)}(v)$. Let $\xi \in N^P_{T^B_S(x)}(v)$, then for every $w\in T^B_S(x)$,
   \[
   \li \xi,w\ri=\li \xi,w-v\ri+\li \xi,v\ri\leq 0.
   \]
   Thus $\xi\in \left(T^B_S(x)\right)^{\circ}$ and by Lemma \ref{pol} it follows that $\xi\in N^P_S(x)$ .

   Hence $\exp_{x}^{-1}y-v\in N^P_S(x)$ and so we have
   \[
   \left\li \frac{\exp_{x}^{-1}y-v}{\norm \exp_{x}^{-1}y-v\norm}, \exp_x^{-1}y\right\ri\leq \varphi(x)d^2(x,y).
   \]
   This implies that $\norm \exp_{x}^{-1}y-v\norm \leq \varphi(x)d^2(x,y)$ which completes the proof. \qed

 \end{proof}

%%%%%%%%%%%%%%%%%%%%%%%%%%%%%%%%%%%%%%%%%%%%%%%%%%%%%%%%%%%%%%%%%%%%%%%%%%%%%%%%%%%%%%%%%%%%%%%%%%%%%%%%%%%%%%%%%%%%%%%%%%%%%%%%%%
%%%%%%%%%%%%%%%%%%%%%%%%%%%%%%%%%%%%%%%%%%%%%%%%%%%%%%%%%%%%%%%%%%%%%%%%%%%%%%%%%%%%%%%%%%%%%%%%%%%%%%%%%%%%%%%%%%%%%%%%%%%%%%%%%%

\re {We are now ready to prove} that the projection map $P_S$ is locally Lipschitz on an open set containing $S$, where $S$ is a prox-regular subset of $M$. In \cite{Hosseini on}, this property of prox-regular sets is verified in the special case in which $M$ is a Hadamard manifold.

 Recall that the Hessian of a $C^2$ function $\psi$ on $M$ is defined by
\[
\hess \psi(x)(v,w):=\left\li \nabla_X \nabla \psi, Y\right\ri(x),
\]
for every $x\in M$ and $v,w\in T_xM$ where $X$, $Y$ are any vector fields such that $X(x)=v$ and $Y(x)=w$ and $\nabla \psi$ denotes the gradient of $\psi$.

\begin{lemma}\label{hess}
Let $M$ be a Riemannian manifold and $x\in M$. Assume that  $R>0$ and $k_0>0$ are given such that $\mid k\mid \leq k_0$ for every sectional curvature $k$ on $B(x,R)$. Then the function $\psi(z):=d^2(x,z)$ is smooth on $B(x,r)$ for every $r>0$ with $r<\min\left\{r(x),R,\frac{\pi}{2\sqrt{k_0}}\right\}$ and
  \[
    \hess \psi(z)(w)^2\geq c(z)\norm w \norm^2,
   \]
   for every $z\in B(x,r)$ and $w\in T_zM$, where
   \[
   c(z)=\min\left\{2,2\sqrt{k_0}\:d(x,z)\cot\left(\sqrt{k_0}\:d(x,z)\right)\right\}.
   \]
\end{lemma}
\begin{proof}
  Let $z\in B(x,r)$ and $w\in T_zM$. Thus according to the proof of \cite[Proposition 2.2]{Fry}, we have
  \[
  \hess \psi(z)(w)^2=2l\left\li D_tX(l),X(l)\right\ri,
  \]
  where $l=d(x,z)$, $X$ is the unique Jacobi field along $\gamma$ with the property that $X(0)=0$ and $X(l)=w$ and $\gamma$ is the unique minimizing geodesic, parameterized by arc length, such that $\gamma(0)=x$ and $\gamma(l)=z$.

  Let $w=w^\top+w^\bot$ be the orthogonal decomposition of $w$  where $w^\top$ is tangent to $\gamma$ and $w^\bot$  is orthogonal to $\dot{\gamma}$ at $z$. Using Propositions 2.3 and 2.4 of Chapter IX of \cite{Lang}, the Jacobi field $X$ can be decomposed into $X=X^\top+X^\bot$ where $X^\top$ and $X^\bot$ are Jacobi fields along $\gamma$ with the property that $X^\top$ and $D_tX^\top$ are tangent to $\gamma$ and $X^\bot$ and $D_tX^\bot$ are orthogonal to $\gamma$. So $X^\top(l)=w^\top$ and $X^\bot(l)=w^\bot$ and using the proof of \cite[Proposition 2.2]{Fry},
  \[
  \begin{array}{ll}
    \hess \psi(z)(w)^2&=   2l\left\li D_tX(l),X(l)\right\ri \\
     & = 2l\left\li D_tX^\top(l),X^\top(l)\right\ri+2l\left\li D_tX^\bot(l),X^\bot(l)\right\ri\\
     & \geq 2l\left(\frac{1}{l}{\norm w^\top\norm}^2\right)+2l\sqrt{k_0}\cot\left(l\sqrt{k_0}\right){\norm w^\bot\norm}^2 \\
     & = 2{\norm w^\top\norm}^2+2l\sqrt{k_0}\cot\left(l\sqrt{k_0}\right){\norm w^\bot\norm}^2 \\
     & \geq c(z){\norm w\norm}^2
    \end{array}
    \]
  where $c(z)=\min\left\{2,2\sqrt{k_0}\:d(x,z)\cot\left(\sqrt{k_0}\:d(x,z)\right)\right\}$. \qed

\end{proof}
\begin{theorem}\label{lip}
  Suppose that $S$ is a closed prox-regular subset of a Riemannian manifold $M$. Then $P_S$ is locally Lipschitz on an open set $V$ containing  $S$.
\end{theorem}
\begin{proof}
\blu{Since $S$ is prox-regular, there exists a continuous function $\varphi : S\rightarrow [0,\infty)$ such that $S$ is   $\varphi$-convex.}
  Let $x\in S$ and $R>0$ be such that $R<r(x)$ and $B(x,R)$ has compact closure and $B(x,R)\subseteq \unp (S)$. Suppose that $k_0>0$ and $\rho >0$ are two constants  such that $\mid k\mid \leq k_0$ for every sectional curvature $k$ on $B(x,R)$ and $\varphi(z)\leq \rho$ for every $z\in B(x,r)\cap S$. Consider $\rbar>0$ given by $\rbar\leq r(z)$ for every $z\in B(x,R)$.

   Let $a\in \mathbb{R}$  be \re{the} solution of the equation $2t\cot (t)=1$ on the interval $\left(0,\frac{\pi}{2}\right)$. So we have $2t\cot (t)>1$ for every $t\in (0,a)$. We now choose $r>0$ such that
   \[
   r<\min\left\{\frac{R}{2}, \rbar, \frac{1}{4\rho}, \frac{a}{\sqrt{k_0}}\right\}.
   \]
   We show that $P_S$ is Lipschitz on $B(x,r)$. To this end, let $x_1,x_2\in B(x,r)$. The \re{case} $x_1,x_2\in S$ is trivial, hence we suppose that $x_1\notin S$. Similar to the proof of \cite[Theorem 3.13]{Hosseini on}, we consider the following optimization problems
   \begin{equation}\label{i}
     \min_{s\in S\cap B(x,R)}d^2\left(x_1,s\right)=\min_{v\in \expxv\left(S\cap B(x,R)\right)}d^2\left(x_1,\expx v\right),
   \end{equation}
   \begin{equation}\label{ii}
     \min_{s\in S\cap B(x,R)}d^2\left(x_2,s\right)=\min_{v\in \expxv\left(S\cap B(x,R)\right)}d^2\left(x_2,\expx v\right).
   \end{equation}
   Let $P_S(x_1)=s_1$ and $P_S(x_2)=s_2$, hence $s_1\in B(x,R)$ and $s_1$ is the optimal solution of (\ref{i}). Moreover, $s_1\in B(x_1,r)\subseteq B(x,R)$ because $x\in S$ and $d(x_1,s_1)\leq d(x_1,x)<r$.

   We claim that there exists a positive constant $\sigma$ such that
   \[
   d^2\left(x_1,s\right)\geq d^2\left(x_1,s_1\right)+\sigma d^2 \left(s,s_1\right),
   \]
   for every $s\in S\cap B(x_1,r)$.

   Let $s\in S\cap B(x_1,r)$ and $\gamma(t)=\exp_{s_1}\left(t\exp^{-1}_{s_1}s\right)$ be the unique geodesic joining $s_1$ and $s$ which is entirely in $B(x_1,r)$.  We now define $\psi(z):= d^2(x_1,z)$ for every $z\in M$. Then using the Taylor expansion, there exists $t_0\in (0,1)$ such that
  \begin{equation}\label{D1}
    d^2(x_1,s)=d^2(x_1,s_1)-2\li \exp_{s_1}^{-1}x_1, \exp_{s_1}^{-1}s\ri + \frac{1}{2}\;\hess \psi(x_0)(v_0)^2,
  \end{equation}
  where $x_0=\gamma(t_0)$ and $v_0=\dot{\gamma}(t_0)$. By Lemma \ref{hess},
  \[
  \hess \psi(x_0)(v_0)^2\geq c(x_0)\norm v_0\norm^2=c(x_0)\;d^2(s,s_1),
  \]
   where $c(x_0)=\min\left\{2,2\sqrt{k_0}\:d(x_1,x_0)\cot\left(\sqrt{k_0}\:d(x_1,x_0)\right)\right\}$. Since $x_0 \in B(x_1,r)$, by the choice of $r$ we have %$\sqrt{k_0}\:d(x_1,x_0)<b$ and so
   \[
   2\sqrt{k_0}\:d(x_1,x_0)\cot\left(\sqrt{k_0}\:d(x_1,x_0)\right)>1,
   \]
    and so $c(x_0)>1$. Moreover, $\exp_{s_1}^{-1}x_1\in N^P_S(s_1)$, hence
    \[
    \begin{array}{ll}
      \li \exp_{s_1}^{-1}x_1, \exp_{s_1}^{-1}s\ri & \leq \varphi(s_1)d(x_1,s_1)d^2(s,s_1) \\
       & \leq \rho rd^2(s,s_1).
    \end{array}
    \]
    Therefore (\ref{D1}) turns into
     \begin{equation}\label{D2}
        d^2(x_1,s)\geq d^2(x_1,s_1)+\left(\frac{1}{2}-2\rho r\right)d^2(s,s_1).
     \end{equation}
    We put $\sigma=\left(\frac{1}{2}-2\rho r\right)$, hence our choice of $r$ guarantees that $\sigma>0$ and the proof of the claim is complete.

     Suppose that $\expx(w_i)=s_i$ for $i=1,2$, then (\ref{D2}) implies that
     \[
     d^2(x_1,\expx(v))\geq d^2(x_1,\expx(w_1))+\frac{\sigma}{{c_1}^2} \norm v-w_1\norm^2,
     \]
     for every $v\in \expx^{-1}\left(S\cap B(x_1,r)\right)$, where $c_1$ is the Lipschitz constant of $\expx^{-1}$ on $B(x,R)$.

     Let $c_2$ be a Lipschitz constant of $\expx$ on $B(0_x,R)$, then by Shapiro's variational principle we finally get
     \[
     \begin{array}{ll}
       d\left(P_S(x_1),P_S(x_2)\right) & = d\left(\expx w_1, \expx w_2\right) \\
        & \leq c_2 \norm w_1-w_2\norm \\
        & \leq \frac{2\kappa c_2{c_1}^2}{\sigma}\;d\left(x_1,x_2\right),
     \end{array}
     \]
     where $\kappa$ is a positive constant such that $2\kappa \;d\left(x_1,x_2\right)$ is a Lipschitz constant of the function $f(v)=d^2(x_1,\expx(v))-d^2(x_2,\expx(v))$ on the neighborhood $W:=\expx^{-1}\left( B(x_1,r)\right)$ of $w_1$. \qed
\end{proof}

We recall that a continuous map $r: X\ra A$ from a topological space $X$ to a subspace $A$ of $X$ is said to be a retraction if the restriction  to $A$ of $r$ is the identity map. A subset $S$ of a Riemannian manifold $M$ is called $\mathcal{L}$-retract if there exist a neighborhood $V$ of $S$, a retraction $r:V\ra S$ and a positive constant $L$ such that
\[
d\left(x,r(x)\right)\leq Ld_S(x), \quad \forall x\in V.
\]

\begin{proposition}
  If $S$ is a prox-regular subset of $M$, then $S$ is $\mathcal{L}$-retract with $L=1$.
\end{proposition}
\begin{proof}
  According to Theorem \ref{lip}, the projection map $P_S:V\ra S$ is a locally Lipschitz retraction from a neighborhood $V$ of $S$ to $S$. \qed
\end{proof}

%%%%%%%%%%%%%%%%%%%%%%%%%%%%%%%%%%%%%%%%%%%%%%%%%%%%%%%%%%%%%%%%%%%%%%%%%%%%%%%%%%%%%%%%%%%%%%%%%%%%%%%%%%%%%%%%%%%%%%%%%%%%%%%%%%%%%%%%%%%%%%%
%%%%%%%%%%%%%%%%%%%%%%%%%%%%%%%%%%%%%%%%%%%%%%%%%%%%%%%%%%%%%%%%%%%%%%%%%%%%%%%%%%%%%%%%%%%%%%%%%%%%%%%%%%%%%%%%%%%%%%%%%%%%%%%%%%%%%%%%%%%%%%%%
\section{Directional differentiability of the metric projection at a boundary point}\label{sec5}

 \re{In this section by applying  Shapiro's variational principle, we investigate the directional differentiability of the projection map $P_S$ at the  boundary points of $S$ where $S$ is  a prox-regular subset of a Riemannian manifold $M$. Let us  recall   the definition of directional differentiability for maps between two Riemannian manifolds.}
\begin{definition}
Let $f:M\ra N$ be a map between two Riemannian manifolds, $x\in M$ and $\left(V,\phi\right)$ be a chart of $N$ at the point $f(x)$.  We define the directional derivative of $f$ at $x$ in the direction $v\in T_xM$ as
\[
f'(x;v):=\lim_{t\ra 0^+}\frac{\phi\left(f\left(\expx(tv)\right)\right)-\phi\left(f(x)\right)}{t},
\]
when the limit exists.

Moreover, the map $f$ is said to be directionally differentiable at $x$ if the directional derivative $f'(x;v)$ exists for all $v\in T_xM$.
\end{definition}

In fact, $f'(x;v)$ is the right-handed derivative of the curve $\gamma(t):=f\left(\expx(tv)\right)$ at $t=0$.

\begin{theorem}\label{ddp}
 Let $S$ be a prox-regular subset of $M$ and $x\in S$.
 %and suppose that for every $w\in T^B_S(x_0)$,
 %\begin{equation}\label{q1}
 %  \lim_{t\ra 0^+}\frac{d\left(\exp_{x_0}(tw), S\right)}{t}=0.
 %\end{equation}
  Then $P_S$ is directionally differentiable at $x$ and for every $v\in T_{x}M$
 \[
 P_S'\left(x;v\right)=P_{T^B_S\left(x\right)}(v),
 \]
 where $P_{T^B_S\left(x\right)}$ denotes the metric projection to $T^B_S\left(x\right)$.
\end{theorem}
\begin{proof}
\blu{Prox-regularity of $S$ implies the existence of a continuous function $\varphi : S\rightarrow [0,\infty)$ such that $S$ is $\varphi$-convex.}
  Let $B\left(x,r\right)\subseteq \unp(S)$ be a convex ball with compact closure  and $v\in T_{x}M$. \re{We are going to show that}
  \[
  \lim_{t\ra 0^+}\frac{\exp_{x}^{-1}\left(P_S\left(\exp_{x}(tv)\right)\right)}{t}=P_{T^B_S\left(x\right)}(v).
  \]
  Since $T^B_S\left(x\right)$ is a closed convex cone in $T_{x}M$, we have
  \[
  P_{T^B_S\left(x\right)}(tv)=tP_{T^B_S\left(x\right)}(v) \quad \forall t\geq 0,
  \]
 and so equivalently we must prove that
  \[
  \lim_{t\ra 0^+}\frac{\norm \exp_{x}^{-1}\left(P_S\left(\exp_{x}(tv)\right)\right)-P_{T^B_S\left(x\right)}(tv)\norm}{t}=0.
  \]
  This means that
  \[
  \norm \exp_{x}^{-1}\left(P_S\left(\exp_{x}(tv)\right)\right)-P_{T^B_S\left(x\right)}(tv)\norm=o(t).
  \]
  To this end, let $t>0$ be given such that $t<\frac{r}{2\norm v\norm}$ and consider the following optimization problems
  \begin{equation}\label{i}
    \min_{w\in T^B_S(x)}\norm w-tv\norm ^2
  \end{equation}
  and
  \begin{equation}\label{ii}
    \min_{y\in S\cap B\left(x,r\right)}d^2\left(y,\exp_{x}tv\right)=\min_{w\in \exp_{x}^{-1}\left( S\cap B\left(x,r\right)\right)}d^2\left(\exp_{x}w,\exp_{x}tv\right).
  \end{equation}
  Note that $\exp_{x}tv\in B\left(x,r\right)\subseteq \unp(S)$, then we get $\xbar=P_S\left(\exp_{x}(tv)\right)$. Since $x\in S$, we have
  \[
  \begin{array}{ll}
    d\left(\xbar,x\right) & \leq d\left(\xbar,\exp_{x}(tv)\right)+d\left(\exp_{x}(tv),x\right) \\
     & \leq 2d\left(\exp_{x}(tv),x\right)=2t\norm v\norm<r,
  \end{array}
  \]
  so $\xbar\in B\left(x,r\right)$. Hence  $\bar{v}=\exp_{x}^{-1}\left(P_S\left(\exp_{x}(tv)\right)\right)$ is the optimal solution of (\ref{ii}).

  Let $v^{*}$ be the optimal solution of (\ref{i}), then $v^{*}=P_{T^B_S\left(x\right)}(tv)$. Furthermore using the proof of \cite[Theorem 3.1]{Shapiro1},
  \[
  \norm w-tv\norm ^2\geq \norm v^*-tv\norm ^2+\norm w-v^*\norm ^2,
  \]
  and (\ref{c1}) is the case for $\alpha=1$. We take $\rbar:=2t\norm v\norm$ and $W:=B(0,\rbar)\subseteq T_{x}M$, then $\bar{v}, v^*\in W$ and by Shapiro's variational principle,
  \[
  \norm \bar{v}-v^* \norm \leq \vartheta(t),
 \]
 where
 \[
 \vartheta(t)=\kappa(t)+2\delta_1(t)+\left(k_1(t)\delta_1(t)+k_2(t)\delta_2(t)\right)^{1/2},
 \]
 and $\kappa(t)$ is a Lipschitz constant of the function
 \[
 h_t(w):=d^2\left(\exp_{x}w,\exp_{x}tv\right)-\norm w-tv\norm ^2,
 \]
 on $W$. Moreover, $k_1(t)$ and $k_2(t)$ are Lipschitz constants of the functions $f_t(w):=\norm w-tv\norm ^2$ and $g_t(w):=d^2\left(\exp_{x}w,\exp_{x}tv\right)$ on $W$, respectively and
 \[
 \begin{array}{ll}
 \delta_1(t)=&\sup\left\{d\left(\exp_{x}^{-1}y, T^B_S(x)\right):y\in S\cap B\left(x,\rbar\right)\right\},\\
 \delta_2(t)=&d\left(v^*,\exp_{x}^{-1}\left( S\cap B\left(x,\rbar\right)\right)\right).
 \end{array}
 \]

 We now \re{show} that $\vartheta(t)=o(t)$. Indeed,
 \[
 k_1(t)\leq 6t\norm v\norm
 \]
 and
 \[
 \begin{array}{ll}
 k_2(t)&=\max_{w\in \overline{W}}\norm \li -2\exp_{\exp_{x}(w)}^{-1}\exp_{x}(tv), d\exp_{x}(w)\ri\norm\\
       &\leq \max_{w\in \overline{W}}\left( 2d\left(\exp_{x}(w), \exp_{x}(tv)\right)\norm d\exp_{x}(w)\norm\right)\\
      &\leq 6lt\norm v\norm,
 \end{array}
 \]
 where $l=\max_{w\in \overline{W}}\norm d\exp_{x}(w)\norm$. Hence $k_1(t)\ra 0$ and $k_2(t)\ra 0$ as $t\ra 0^+$.
 Moreover, by Theorem \ref{o}, for every $y\in S\cap B(x,r)$
 \[
 d\left(\exp_{x}^{-1}y, T^B_S(x)\right)\leq \varphi(x)d^2\left(x,y\right).
 \]
 So we have
 \[
 \lim_{\stackrel{y\ra x}{y\in S}}\frac{d\left(\exp_{x}^{-1}y, T^B_S(x)\right)}{d\left(x,y\right)}=0,
 \]
 and this implies that $\delta_1(t)=o(t)$. Also, %it follows from (\ref{q1}) that
  $\delta_2(t)=o(t)$  since
  \[
 \begin{array}{ll}
 \delta_2(t)=& d\left(tv_0,\exp_{x}^{-1}\left( S\cap B\left(x,\rbar\right)\right)\right)
 \leq c_1d\left(\exp_{x}(tv_0),S\right),
 \end{array}
 \]
  where $v_0=P_{T^B_S(x)}(v)$ and $c_1$ is a Lipschitz constant of $\exp_{x}^{-1}$ on $B\left(x,\rbar\right)$.

  It remains only to verify that $\kappa(t)=o(t)$. Indeed, for every $w\in B(0,r)$ and $z\in T_{x}M$,
  \[
  \nabla h_t(w)(z)=-2\left \li \exp_{\exp_{x}(w)}^{-1}\exp_{x}(tv), d\exp_{x}(w)z\right\ri-2\li w-tv,z\ri.
  \]
  For fixed $w,z$ we define
  \[
  F(t)= \left\li \exp_{\exp_{x}(w)}^{-1}\exp_{x}(tv), d\exp_{x}(w)z\right\ri,
  \]
  for every $t$ with $|t|<\frac{r}{2\norm v\norm}$. The Taylor expansion gives
  \[
  F(t)=F(0)+F'(0)t+o(t) \quad \forall t.
  \]
  The values $F(0)$ and $F'(0)$ is obtained as follows: according to \cite[Lemma 3.5, p. 250]{Lang} we have
  \[
  \begin{array}{ll}
    F(0) & =\left\li \exp_{\exp_{x}(w)}^{-1}x, d\exp_{x}(w)z\right\ri \\
     & = \left\li d\exp_{\exp_{x}(w)}\left(-\dot{\gamma}(1)\right)\left(\exp_{\exp_{x}(w)}^{-1}x\right), z\right\ri,
  \end{array}
  \]
  where $\gamma$ is the geodesic $\gamma(t)=\exp_{x}(tw)$ and hence
  \[
  \dot{\gamma}(1)=-\exp_{\exp_{x}(w)}^{-1}x.
  \]
  For simplicity, let us write $\bar{w}=\exp_{\exp_{x}(w)}^{-1}x$. Thus using \cite[Theorem 3.1]{Lang},
  \[
  d\exp_{\exp_{x}(w)}(\bar{w})(\bar{w})=J(1),
  \]
  where $J$ is the Jacobi field along the geodesic $\beta$ joining $\exp_{x}(w), x$ satisfying the properties $\dot{\beta}(0)=\bar{w}$, $J(0)=0$ and $D_t J(0)=\bar{w}$. In fact,  $\beta (t)=\gamma (1-t)$ and $J(t)=t\dot{\beta}(t)$ and so $J(1)=-w$ and $F(0)=-\li w,z\ri$.

  Also we have
   \[
  \begin{array}{ll}
  F'(0)&=\left\li d\exp_{\exp_{x}(w)}^{-1}(x)v,d\exp_{x}(w)z\right\ri\\
  &= \left\li d\exp_{\exp_{x}(w)}\left(\exp_{\exp_{x}(w)}^{-1}x\right)\left(d\exp_{\exp_{x}(w)}^{-1}(x)v\right),z\right\ri\\
  &= \left\li d\left(\exp_{\exp_{x}(w)}o\exp_{\exp_{x}(w)}^{-1}\right)(x)v,z\right\ri\\
  &= \li v,z\ri.
  \end{array}
  \]
  It follows that $ \nabla h_t(w)(z)=o(t)$ for every $w\in B(0,r)$ and $z\in T_{x}M$. This implies that $\kappa(t)=o(t)$. \qed
\end{proof}

Using Theorem \ref{ddp}, we  obtain the following characterization of Bouligand tangent cone to a prox-regular set.

\begin{corollary}\label{btc}
  Let $S$ be a closed prox-regular subset of $M$ and $x\in S$. Then $v\in T^B_S(x)$ if and only if there exists a continuous curve $\alpha:[0,\varepsilon)\ra S$ such that $\alpha(0)=x$ and $\dot{\alpha}(0^+)=v$, where $\dot{\alpha}(0^+)$ is the right-handed derivative of $\alpha$ at $0$.
\end{corollary}
\begin{proof}
  Let $v\in T^B_S(x)$. We choose $\varepsilon>0$ such that $\exp_x(tv)\in \unp(S)$ for all $t\in [0,\varepsilon)$.  We now define
  \[
  \alpha(t):=P_S\left(\exp_x(tv)\right)\qquad \forall t\in [0,\varepsilon).
  \]
  Then by Theorem \ref{ddp},
  \[
  \dot{\alpha}(0^+)= P_S'\left(x;v\right)=P_{T^B_S\left(x\right)}(v)=v.
  \]
  The proof of the converse statement is straightforward. \qed
\end{proof}

%\begin{example}
%Suppose that $S$ is a $\varphi$-convex subset of $M$ such that $\partial S$ is a  $C^2$ submanifold of $M$. If $\partial S$ is totally geodesic, then using this property and Proposition \ref{tbss}, for every $x\in \partial S$ and $v\in T^B_S(x)$ there exists $\varepsilon=\varepsilon(x,v)>0$ such that $\exp_x(tv)\in S$ for all $t\in [0,\varepsilon)$. Therefore the condition (\ref{q1}) holds and $P_S$ is directionally differentiable on $\partial S$.
%\end{example}
%%%%%%%%%%%%%%%%%%%%%%%%%%%%%%%%%%%%%%%%%%%%%%%%%%%%%%%%%%%%%%%%%%%%%%%%%%%%%%%%%%%%%%%%%%%%%%%%%%%%%%%%%%%%%%%%%%%%%%%%%%%%%%%%%%
%%%%%%%%%%%%%%%%%%%%%%%%%%%%%%%%%%%%%%%%%%%%%%%%%%%%%%%%%%%%%%%%%%%%%%%%%%%%%%%%%%%%%%%%%%%%%%%%%%%%%%%%%%%%%%%%%%%%%%%%%%%%%%%%%%
\section{Minimizing curves in prox-regular sets}\label{sec6}
 Our goal in this section is to derive a necessary condition for a curve $\gamma$ to be a minimizing curve between its endpoints in a prox-regular set. %whose boundary is a $C^2$ submanifold of $M$.
  To this end, we employ the first variation formula.
 Let \re{$S\subseteq M$}  be a closed prox-regular set whose boundary is a $C^2$ Riemannian submanifold of $M$.

  In this situation, a continuous map $\gamma:[a,b]\ra M$ is called a piecewise regular curve if it is a piecewise $C^2$ curve with nonzero derivatives. Moreover,
 by an admissible curve we mean a piecewise regular curve $\gamma:[a,b]\ra M$ which is entirely in $S$. An admissible curve $\gamma$ in $S$ is said to be minimizing if $\mathcal{L}\left(\gamma\right)\leq  \mathcal{L}\left(\tilde{\gamma}\right)$ for  all admissible curves $\tilde{\gamma}$   with the same endpoints where $\mathcal{L}\left(\gamma\right)$ denotes the length of $\gamma$ in $M$.

 An admissible family of curves in $S$ is a continuous map $\Gamma: [0,\varepsilon)\times [a,b]\ra M$ with the property that $\Gamma(s,t)\in S$ for all $(s,t)\in [0,\varepsilon)\times [a,b]$ and there exists a partition $a=a_0<\cdots<a_k=b$ of $[a,b]$ such that $\Gamma|_{[0,\varepsilon)\times [a_{i-1},a_i]}$ is $C^2$ for every $i=1,\ldots, k$. A variation of an admissible curve $\gamma : [a, b] \ra M$ is an admissible
family $\Gamma$ in $S$ such that $\Gamma(0,t) = \gamma(t)$ for all $t \in [a, b]$ and if in addition $\Gamma(s,a) = \gamma(a)$ and $\Gamma(s,b) = \gamma(b)$  for all $s\in [0,\varepsilon)$, then it is called a
proper variation.

Recall that if $\Gamma$ is a variation of $\gamma$, then the piecewise $C^1$ vector field $V$ along $\gamma$ defined by $V (t) = \frac{d}{ds}|_{s=0^{+}} \Gamma(s, t)$ is called the variation field of $\Gamma$, where $\frac{d}{ds}|_{s=0^{+}}$ denotes the right-handed derivative of $\Gamma(.,t):[0,\varepsilon)\ra M$ at $s=0$. Note that according to Corollary \ref{btc}, if $V$ is the variation field of a variation along $\gamma$, then
\[
V(t)\in T^B_S\left(\gamma(t)\right)\quad \forall t\in [a,b].
\]

In the following, we investigate when a vector field along an admissible curve $\gamma$ is the variation field of a variation of $\gamma$.

 \begin{lemma}\label{tbss}
    Suppose that the closed set $S$ is prox-regular and $\partial S$ is a $C^2$ submanifold of $M$. If $x\in \partial S$ and $v\in \left(T^B_S(x)\setminus T_x\partial S\right)\cup \{0\}$, then there exists $\varepsilon>0$ such that $\exp_x(tv)\in S$ for all $t\in [0,\varepsilon)$.

 \end{lemma}
 \begin{proof}
    Assuming the contrary,  there exists a sequence $\{t_n\}$ such that $t_n\downarrow 0$ and $\exp_x\left(t_nv\right)\in S^c\subseteq \hat{S}$.
     Moreover,
     \[
     v=\lim_{n\ra\infty}\frac{\exp_x^{-1}\left(\exp_x\left(t_nv\right)\right)}{t_n}.
      \]
      Hence $v\in T^B_{\hat{S}}(x)$ and so by Theorem \ref{tbsh}, we have $v\in T_x\partial S$. This contradiction completes the proof. \qed
 \end{proof}

 \begin{lemma}\label{vari}
    Suppose that $\gamma: [a,b]\ra M$ is an admissible curve and $V$ is a piecewise $C^2$ vector field along $\gamma$. If for any $t\in [a,b]$ with  $\gamma(t)\in \partial S$ we have
   \[
   V(t)\in \left(T^B_S\left(\gamma(t)\right)\setminus T_{\gamma(t)}\partial S\right)\cup \{0\},
    \]
    then $V$ is the variation field of a variation $\Gamma$ of $\gamma$.
 \end{lemma}
 \begin{proof}
      Lemma \ref{tbss} along with the compactness of $[a,b]$ imply that there exists $\varepsilon>0$ such that the map $\Gamma: [0,\varepsilon)\times [a,b]\ra M$ defined by $\Gamma(s,t):= \exp_{\gamma(t)}\left(sV(t)\right)$ is the desired variation of $\gamma$ in $S$. \qed
 \end{proof}

    The following theorem gives a necessary condition for a curve to be minimizing in $S$.

\begin{theorem}
Let $\gamma:[a,b]\ra M$ be a unit speed admissible curve. If $\gamma$ is  minimizing in $S$, then
\begin{equation}\label{nece}
D_t\dot{\gamma}(t)\in N^P_S\left(\gamma(t)\right),
 \end{equation}
for every $t\in [a,b]$ except for finitely many points.
 \end{theorem}
 \begin{proof}
 If the interior of $S$ is empty, then (\ref{nece}) evidently holds, since in this case $N^P_S(x)= T^{\bot}_{x}\partial S$ for every $x\in S$. Therefore we assume that the interior of $S$ is nonempty.

 Let $a=a_0<a_1<\cdots<a_k=b$ be a partition of $[a,b]$ such that $\gamma$ is $C^2$ on each subinterval $[a_{i-1},a_i]$ and $t_0\in [a,b]$ be such that $t_0\neq a_i$ for each $i$. We suppose  that $t_0\in (a_{j-1},a_j)$ for some $j$, $1\leq j\leq k$ and we get $x_0:=\gamma\left(t_0\right)$. If $x_0\in S^{\circ}$, then there exist an open neighborhood $U$ of $x_0$ in $M$ and a positive number $\delta$ such that $\gamma(t)\in U\subseteq S$ for all $t\in I_0:=[t_0-\delta, t_0+\delta]$. Hence $\gamma|_{I_0}$ is  minimizing in $M$ and this implies that $D_t\dot{\gamma}\left(t_0\right)=0$.

 Assume that $x_0\in \partial S$ and let the open neighborhood $U$ and the submersion $\psi: U\ra \mathbb{R}$ be the ones applied in the proof of Theorem \ref{C2}.   Clearly $D_t\dot{\gamma}\left(t_0\right)=0$ or there is a positive number $\epsilon$ such that $\gamma(t)\in U\cap\partial S$ for all $t\in I:=[t_0-\epsilon, t_0+\epsilon]\subset (a_{j-1},a_j)$. So it suffices to check that (\ref{nece}) holds in the latter case.
       Indeed, $\gamma|_{I}$ is minimizing in the Riemannian submanifold $\partial S$ of $M$. Then we have
 \begin{equation}\label{tbot}
 D_t\dot{\gamma}\left(t\right)\in T^{\bot}_{\gamma(t)}\partial S\qquad \forall t\in I.
 \end{equation}
 If $U\cap S^{\circ}=\emptyset$ (or there exists a neighborhood $V\subseteq U$ of $x$ such that $V\cap S^{\circ}=\emptyset$), then $N^P_S(x)= T^{\bot}_{x}\partial S$. Otherwise, in order to deduce that $ D_t\dot{\gamma}\left(t_0\right)\in N^P_S\left(\gamma\left(t_0\right)\right)$, by Lemma \ref{pol} it suffices to show that
 \[
\left\li D_t\dot{\gamma}\left(t_0\right), v\right\ri \leq 0\qquad \forall v\in T^B_S\left(\gamma\left(t_0\right)\right).
 \]
If this fails to hold, then there is $\eta\in T^B_S\left(\gamma\left(t_0\right)\right)$ such that
\[
\left\li D_t\dot{\gamma}\left(t_0\right), \eta\right\ri>0.
\]
Thus the inclusion (\ref{tbot}) implies that $\eta\notin T_{\gamma\left(t_0\right)}\partial S$.

We now construct a vector field along $\gamma$ such that for any $t\in [a,b]$ with $\gamma(t)\in \partial S$,
\[
V(t)\in \left(T^B_S\left(\gamma(t)\right)\setminus T_{\gamma(t)}\partial S\right)\cup \{0\}.
\]
We define $\overline{V}(t):=L^{\gamma}_{t_0t}\eta$ and $g(t):=\left\li \overline{V}(t), \nabla\psi\left(\gamma(t)\right)\right\ri$  for all $t\in I$. According to Theorem \ref{C2}, $N^P_S\left(\gamma(t)\right)= \cone \left\{\nabla\psi\left(\gamma(t)\right)\right\}$ for all $t\in I$. Then $g\left(t_0\right)<0$ and  the continuity of $g$ implies that  $g(t)<0$ on a possibly smaller neighborhood of $t_0$. It follows that
\[
\overline{V}(t)\in T^B_S\left(\gamma(t)\right)\setminus T_{\gamma(t)}\partial S \qquad \forall t\in I,
\]
without loss of generality.
By shrinking $I$ if necessary, we can assume that
\begin{equation}\label{vdt}
 \left\li \overline{V}(t), D_t\dot{\gamma}\left(t\right)\right\ri>0\qquad \forall t\in I.
\end{equation}
We choose a bump function $\phi \in C^{\infty}(\mathbb{R})$ with support in $I$ such that $\phi(t)\equiv1$ on $[c_1,c_2]$, where $c_1,c_2$ is such that $t_0-\epsilon<c_1<t_0<c_2<t_0+\epsilon$. We now define $V(t):=\phi(t)\overline{V}(t)$ for all $t\in [a,b]$. Then $V$ is the desired vector field along $\gamma$.

Applying Lemma \ref{vari} to the vector field $V$ along $\gamma$ gives rise to  a variation $\Gamma:[0,\varepsilon)\times [a,b]\ra M$ of $\gamma$ in $S$ such that $V$ is its variation field. Since $\gamma$ is minimizing in $S$, for all $s\in [0,\varepsilon)$ we have $\mathcal{L}\left(\Gamma_s\right)\geq \mathcal{L}\left(\Gamma_0\right)$  where $\Gamma_s$ is an admissible curve on $[a,b]$ defined by $\Gamma_s(t):=\Gamma(s,t)$. This implies that $\frac{d}{ds}|_{s=0^{+}}\mathcal{L}\left(\Gamma_s\right)\geq 0$. \re{Using the first variation formula (see \cite[Theorem 6.3]{lee})} we conclude that
\[
\int_{a}^{b} \left\li V(t), D_t\dot{\gamma}(t)\right\ri dt+\sum_{i=1}^{k-1}\left\li V(a_i), \Delta_i\dot{\gamma}\right\ri\leq 0,
\]
where $\Delta_i\dot{\gamma}:=\dot{\gamma}\left(a_i^+\right)-\dot{\gamma}\left(a_i^-\right)$. Since $V(a_i)=0$ for each $i=1,\ldots,k$, we have
\[
\int_{a}^{b} \left\li V(t), D_t\dot{\gamma}(t)\right\ri dt\leq 0.
\]
On the other hand,
\[
\int_{a}^{b} \left\li V(t), D_t\dot{\gamma}(t)\right\ri dt\geq \int_{c_1}^{c_2} \left\li \overline{V}(t), D_t\dot{\gamma}(t)\right\ri dt>0,
\]
a contradiction which establishes that
\[
 D_t\dot{\gamma}\left(t_0\right)\in N^P_S\left(\gamma\left(t_0\right)\right).
\]
\qed
\end{proof}
\begin{example}
Let $S^2$ be the $2$-sphere of radius one in $\mathbb{R}^3$ with the round metric $g^\circ$ which is induced from the Euclidean metric on $\mathbb{R}^3$. Consider spherical coordinates $\left(\theta,\phi\right)$ on the subset $S^2-\left\{(x,y,z): x\leq 0,y=0\right\}$ of the sphere defined by
\[
(x,y,z)=\left(\sin \theta \cos \phi, \sin \theta\sin \phi,\cos \theta\right) \quad 0<\theta<\pi,-\pi<\phi<\pi.
\]
It is known that the round metric is $g^\circ=d\theta ^2+\sin ^2 \theta\:d\phi ^2$ in spherical coordinates. Also, Christoffel symbols of $g^\circ$ in spherical coordinates are
\[
\Gamma^{\theta}_{ij}=\left(
\begin{array}{cc}
  0 & 0 \\
  0 & -\sin \theta\cos\theta
\end{array}\right),\qquad
\Gamma^{\phi}_{ij}=\left(
\begin{array}{cc}
  0 & \frac{\cos\theta}{\sin \theta} \\
  \frac{\cos\theta}{\sin \theta} & 0
\end{array}\right).
\]

Let $S$ be the closed subset of $S^2$ which is obtained by removing the sector $\theta_0< \theta \leq\pi$ from $S^2$ where $\frac{\pi}{2}<\theta_0<\pi$.  According to \cite[Theorem 4.18]{convex}, $S$ is a prox-regular subset of $S^2$. %for some continuous function $\varphi:S\ra [0,+\infty)$.

Note that for every $p\in \partial S$, the map $\psi:U\ra \mathbb{R}$ defined by $\psi\left(\theta,\phi\right)=\theta-\theta_0$ is the desired submersion which is used in Theorem \ref{C2}. Hence applying Theorem \ref{C2}, we obtain
  \[
  N^P_S(p)=\cone \{(1,0)\}=\{\lambda\:\partial/\partial \theta:\lambda\geq 0\}.
  \]
Clearly, the unit speed curve $\gamma$ defined by
\[
\gamma(t):=\left(\theta_0,\frac{t}{\sin \theta_0}\right) \qquad \forall t\in I:=\left[-\pi/2\sin \theta_0,\pi/2\sin \theta_0\right],
\]
 is a minimizing curve  in $S$ joining $\left(\theta_0,-\pi/2\right)$ and $\left(\theta_0,\pi/2\right)$. It can be found that the curve $\gamma$ satisfies the necessary condition (\ref{nece}). Indeed, we have
 \[
  D_t\dot{\gamma}\left(t\right)=\left(-\frac{\cos \theta_0}{\sin \theta_0}\right)\frac{\partial}{\partial \theta} \qquad \forall t\in I.
 \]
 Then putting $\lambda:=-\cos \theta_0/\sin \theta_0$, we observe  that $\lambda>0$ and (\ref{nece}) holds.

 On the other hand,  consider another admissible curve $\alpha$ in $S$ joining $\left(\theta_0,-\pi/2\right)$ and $\left(\theta_0,\pi/2\right)$ defined by
 \[
 \alpha(t):=\left\{
 \begin{array}{lc}
   \left(\theta_0-t,-\pi/2\right) & 0\leq t\leq \theta_0-\pi/2 \\
   \left(\pi/2, t-\theta_0\right) & \theta_0-\pi/2\leq t\leq \theta_0+\pi/2 \\
   \left(t-\theta_0,\pi/2\right) & \theta_0+\pi/2\leq t\leq 2\theta_0.
 \end{array}
 \right.
 \]
 Note that $D_t\dot{\alpha}\left(t\right)=0$ and $\alpha$ satisfies (\ref{nece}), but it is not a minimizing curve in $S$, since $\mathcal{L}\left(\gamma\right)=\pi\sin \theta_0<\pi<\mathcal{L}\left(\alpha\right)$.

\end{example}
\begin{example}
  Let $H^2$ be the hyperbolic plane; that is, the upper half-plane in $\mathbb{R}^2$
with the metric $g_H=\left(dx^2+dy^2\right)/y^2$. The Riemannian distance between two points $z_1=\left(x_1,y_1\right), z_2=\left(x_2,y_2\right)$ of $H^2$ is as
\[
d\left(z_1,z_2\right)= 2\ln \frac{\sqrt{\left(x_2-x_1\right)^2+\left(y_2-y_1\right)^2}+\sqrt{\left(x_2-x_1\right)^2+
\left(y_2+y_1\right)^2}}{2\sqrt{y_1y_2}}.
\]

We consider the subset $S=\left\{(x,y):1\leq y\leq 2\right\}$ of $H^2$. It is evident that $S$ is not convex in the Hadamard manifold $\left(H^2,g_H\right)$. The metric projection $P_S$ is obtained as follows:
\[
P_S(x,y)=\left\{
\begin{array}{ll}
  (x,1) & {\textrm{if}}\ \ 0<y<1 \\
  (x,2) & {\textrm{if}}\ \  y>2.
\end{array}
\right.
\]
Since for instance in the case when $y>2$, for every $s\neq x$ we have
%\begin{align*}
\[
  d\left((x,y),(s,2)\right) %& =2\ln \frac{\sqrt{\left(s-x\right)^2+\left(2-y\right)^2}+\sqrt{\left(s-x\right)^2+\left(2+y\right)^2}}{2\sqrt{2y}}\\&
  > \ln \frac{y}{2}=d\left((x,y),(x,2)\right).
 \]
%\end{align*}
Then by \cite[Corollary 4.20]{convex}, $S$ is a prox-regular subset of $H^2$.

Note that for every $z=(s,2)\in \partial S$, the map $\psi:H^2\ra \mathbb{R}$ defined by $\psi\left(x,y\right)=y-2$ is the desired submersion which is needed   in Theorem \ref{C2}. Hence
  \[
  N^P_S(z)=\cone \{4\partial/\partial y\}=\{\lambda\:\partial/\partial y:\lambda\geq 0\}.
  \]
We now consider the unit speed curve $\gamma$ in $S$ defined by
\[
\gamma(t):=\left(2t,2\right) \qquad \forall t\in \mathbb{R}.
\]
So
 $D_t\dot{\gamma}\left(t\right)=2\frac{\partial}{\partial y}$ for all $t\in \mathbb{R}$ and it follows that the curve $\gamma$ satisfies the necessary condition (\ref{nece}).
 \end{example}
%%%%%%%%%%%%%%%%%%%%%%%%%%%%%%%%%%%%%%%%%%%%%%%%%%%%%%%%%%%%%%%%%%%%%%%%%%%%%%%%%%%%%%%%%%%%%%%%%%%%%%%%%%%%%%%%%%%%%%%%%%%%%%%%%%
%%%%%%%%%%%%%%%%%%%%%%%%%%%%%%%%%%%%%%%%%%%%%%%%%%%%%%%%%%%%%%%%%%%%%%%%%%%%%%%%%%%%%%%%%%%%%%%%%%%%%%%%%%%%%%%%%%%%%%%%%%%%%%%%%%

%\begin{acknowledgements}
%The second-named author was supported by a grant from the University of Isfahan.
%\end{acknowledgements}

%%%%%%%%%%%%%%%%%%%%%%%%%%%%%%%%%%%%%%%%%%%%%%%%%%%%%%%%%%%%%%%%%%%%%%%%%%%%%%%%%%%%%%%%%%%%%%%%%%%%%%%%%%%%%%%%%%%%%%%%%%%%%%%%%%%%


\begin{thebibliography}{}

\bibitem{Fry}  Azagra, D.,  Fry, R.: A second order smooth variational principle on Riemannian manifolds. Canad. J. Math. 62,  241-260 (2010)


\bibitem{EFP}  Bangert, V.:  Sets with positive reach. Arch. Math. 38, 54-57 (1982)

\bibitem{Barani}   Barani, A.,  Hosseini, S.,   Pouryayevali, M.R.: On the metric projection onto $\varphi$-convex subsets of Hadamard manifolds. Rev. Mat. Complut. 26, 815-826 (2013)

%\bibitem{Canino1} A. Canino,  Local properties of geodesics on $p$-convex sets, Annali di Matematica Pura ed Applicata, 159.1 (1991), 17-44.
\bibitem{canino1} \blu{Canino, A.:  Local properties of geodesics on $p$-convex sets. Ann. Mat. Pura Appl. 159(1), 17–44 (1991)}

\bibitem{canino}  Canino, A.: On $p$-convex sets and geodesics. J. Differential Equations.  75, 118-157 (1988)

\bibitem{Clark book a}  Clarke, F.H.,  Ledyaev, Yu.S.,  Stern, R.J.,  Wolenski, P.R.: Nonsmooth Analysis and Control Theory. Graduate Texts in Mathematics 178, Springer, New York (1998)

\bibitem{colombo1}  Colombo, G.,  Thibault, L.: Prox-regular sets and applications.  Handbook of Nonconvex Analysis and Applications, D.Y. Gao and D. Motreanu Eds., International Press, Boston, 99-182 (2010)

\bibitem{colombo2}  Degiovanni, M.,  Marino, A.,  Tosques, M.: General properties of $(p,q)$-convex functions and $(p, q)$-monotone operators.  Ricerche Mat. 32, 285-319 (1983)

\bibitem{Do Carmo} do Carmo,  M.P.: Riemannian Geometry. Birkh\"auser, Boston (1992)

%\bibitem{Domokos} A. Domokos,  J. M. Ingram, and M. M. Marsh, Projections onto closed convex sets in Hilbert spaces, Acta Mathematica Hungarica, 152.1 (2017), 114-129.

\bibitem{Federer}  Federer, H.: Curvature measure. Trans. Amer. Math. Soc. 93, 418-491 (1959)


%\bibitem{Fitzpatrick} S. Fitzpatrick, R. R. Phelps, Differentiability of the metric projection in Hilbert space, Transactions of the American Mathematical Society, 270.2 (1982), 483-501.

\bibitem{Gilberg}  Gilbarg, D.,  Trudinger, N.S.: Elliptic Partial Differential Equations of
Second Order. Springer-Verlag, Berlin, New York (1998)

\bibitem{Greene}  Greene, R.E.,  Shiohama, K.: Convex functions on complete noncompact manifolds: topological
structure. Invent. Math. 63, 129-157 (1981)

%\bibitem{Grognet} S. Grognet, Th\'{e}or\`{e}me de Motzkin en courbure n\'{e}gative, Geom. Dedicata 79 (2000), 219-227.

%\bibitem{Holmes} R. B. Holmes,  Smoothness of certain metric projections on Hilbert space, Transactions of the American Mathematical Society, 184 (1973), 87-100.

%\bibitem{Hosseini 1} S. Hosseini,  M. R. Pouryayevali, Generalized gradients and characterization of epi-Lipschitz sets in Riemannian manifolds, Nonlinear Anal. 74 (2011) 3884-3895.

\bibitem{Hosseini on}  Hosseini, S.,  Pouryayevali,  M.R.: On the metric projection onto prox-regular subsets of Riemannian manifolds. Proc. Amer. Math. Soc. 141,  233-244 (2013)

\bibitem{Klein2} Kleinjohann,  N.: Convexity and the unique footpoint property in Riemannian geometry. Arch. Math. 35, 574-582 (1980)

\bibitem{kruskal}  Kruskal, J.: Two convex counterexamples: a discontinuous envelope function and a
non-differentiable nearest-point mapping. Proc. Amer. Math. Soc. 23, 697-703 (1969)

\bibitem{Lang} Lang,  S.: Fundamentals of Differential Geometry.  Graduate Texts in Mathematics 191, Springer-Verlag, New York (1999)

\bibitem{lee}   Lee, J.M.:  Introduction to Riemannian Manifolds.
Graduate Texts in Mathematics 176, Springer, New York (2018)

\bibitem{Leobacher}  Leobacher, G.,   Steinicke,  A.: Existence, uniqueness and regularity of the projection onto differentiable manifolds. arXiv preprint arXiv:1811.10578 (2018)

%\bibitem{Malanowski} K. Malanowski,  Remarks on differentiability of metric projections onto cones of nonnegative functions, Journal of Convex Analysis, 10.1 (2003), 285-294.

\bibitem{RO}  Poliquin, R.A.,  Rockafellar, R.T.: Prox-regular functions in variational analysis. Trans. Amer. Math. Soc.  348, 1805-1838 (1996)

%\bibitem{RO2} R. A. Poliquin, R. T. Rockafellar and L. Thibault, Local differentiability of distance functions, Trans. Amer. Math. Soc.  352 (2000) 5231-5249.

\bibitem{convex}  Pouryayevali, M. R.,  Radmanesh, H.: Sets with the unique footpoint property and $\varphi$-convex subsets of Riemannian manifolds. J. Convex Anal. 26, 617-633 (2019)

%\bibitem{Rataj} J. Rataj, L. Zajicek, On the structure of sets with psitive reach, Math. Nachr., Math. Nachr. (2017) 290 1806–1829.

\bibitem{sakai}  Sakai, T.:  Riemannian Geometry. Translations of Mathematical Monographs 149. American Mathematical Society (1996)

%\bibitem{Shapiro2} A. Shapiro,  Differentiability properties of metric projections onto convex sets, Journal of Optimization Theory and Applications, 169(3) (2016), 953-964.

\bibitem{Shapiro1}  Shapiro, A.: Existence and differentiability of metric projections in Hilbert spaces. SIAM J. Optim. 4(1),  130-141 (1994)

%\bibitem{Thale} C. Th\"ale, 50 years sets with positive reach-A survey, Surveys in Mathematics and its Applications 3  123-165 (2008).

%\bibitem{Villani} C. Villani, \emph{Optimal transport: old and new}, Vol. 338, Springer Science \& Business Media, 2008.

\bibitem{Walter}  Walter, R.: On the metric projection onto convex sets in Riemannian spaces. Arch. Math. (Basel), 25, 91-98 (1974)

\end{thebibliography}
\end{document}